\documentclass{article}

\usepackage{amscd}
\usepackage{amssymb}

\newtheorem{theorem}{Theorem}[section]

\newenvironment{definition}{\medskip 
\noindent {\bf Definition.}}{\mbox{}}

\begin{document}

\centerline{\Large\bf Equivariant Cohomology and Wall Crossing}
\centerline{\Large\bf  Formulas in 
Seiberg-Witten Theory}\bigskip
\centerline{Huai-Dong Cao \& Jian Zhou}


\footnotetext[2]{The authors are supported in part by NSF grant
\#DMS-9632028. }

\begin{abstract}
 We use localization formulas 
in the theory of equivariant cohomology to 
rederive the wall crossing formulas of Li-Liu \cite{Li-Liu} and
Okonek-Teleman \cite{Oko-Tel} for 
Seiberg-Witten invariants.
\end{abstract}

One of the  difficulties in the study of  Donaldson invariants or 
Seiberg-Witten invariants for closed oriented $4$-manifold with 
$b_2^+ = 1$ is that one has to deal with  reducible solutions. There have 
been a lot of work in this direction in the Donaldson theory context (see 
G\"{o}ttsche \cite{Got} and the references therein).
In the Seiberg-Witten theory, the $b_1 = 0$ case were discussed 
in Witten \cite{Wit} and Kronheimer-Mrowka \cite{Kro-Mro}. The general case 
was solved by Li-Liu \cite{Li-Liu}. Very recently, Okonek-Teleman \cite{Oko-Tel}
extended the definition of Seiberg-Witten invariants when 
$b_1 \neq 0$ and obtained a universal wall crossing formula for
the invariants. A common feature in such works is that 
the equations used to define the invariants depend on some parameters. The 
parameter spaces are divided into chambers by walls, where reducible 
solutions can occur. Within the same  chamber,
the invariants do not change. When the parameter changes smoothly from 
 one chamber to  another, the usual 
approach is to examine what happens when one cross the wall. The
result is expressed as a wall crossing formula.

The above complications actually all come from one source: the configuration
spaces are singular. They are 
quotients of contractible spaces by gauge groups, but the 
reducible solutions
and irreducible solutions have different orbit types. This leads one 
to consider
other cohomology theories for the configuration spaces.
For example, Li-Tian \cite{Li-Tia} have used intersection 
homology 
to study the wall crossing phenomenon addressed in Li-Liu \cite{Li-Liu}. 
On the other
hand, there have been many papers in physics literature defining topological
quantum field theories by equivariant cohomology related to  the action of
gauge groups. Such cohomology theories are infinite dimensional in nature,
since the gauge groups are infinite dimensional. 
 In this paper, we use an essentially finite dimensional approach. 
 A well-known procedure is to break the action of the gauge group
 into a free action by an infinite dimensional group, then followed 
 by a finite dimensional compact group. 
 Localization formulas in equivariant
 cohomology theory in the finite dimensional case can then be applied 
 to study the wall crossing.  
 
 The rest of the paper is arranged as follows. $\S 1$ reviews the 
 definition of equivariant cohomology. A localization formula due to
 Kalkman \cite{Kal} and two special cases are discussed in $\S 2$. 
 $\S 3$ and $\S 4$ describe, respectively, how to use the localization formulas to 
 derive the wall crossing formulas for Seiberg-Witten invariants 
 due to Li-Liu \cite{Li-Liu} and Okonek-Teleman \cite{Oko-Tel}.
 
 \section{Equivariant cohomology}

  For simplicity of the presentation, we review only what
 we will use later about equivariant cohomology. For the
 general theory on equivariant cohomology, the reader
 is referred to  \cite{Qui, Ati-Bot, Ber-Get-Ver}.
 We shall only consider the case of an $S^1$-action on
 a compact smooth manifold $W$, with fixed point set $F$. We allow $W$  to
 have boundary, but require that $F \cap \partial W = \emptyset$. 
 The action of $S^1$ generates a vector field $X$ on $W$. In fact, for any
 $x \in W$, if we let $c(t) = \exp (\sqrt{-1} t) \cdot x$ then  $X(x)$ is the tangent
 vector to $c(t)$ at $t = 0$. Denote by $\Omega^*(W)^{S^1}$ the space of 
 differential forms on $W$ fixed under the $S^1$-action.  Let $u$ be an indeterminate 
 of degree $2$ and consider the space $\Omega^*(W)^{S^1} \otimes {\mathbb R}[u]$. 
 Define 
 $$ d_{S^1} = d - u \cdot i_{X} : \Omega^*(W)^{S^1} \otimes {\mathbb R}[u]
 \rightarrow \Omega^*(W)^{S^1} \otimes {\mathbb R}[u]$$ 
 as a derivation, whose action on $u$ is zero and 
 $d_{S^1} \alpha = d \alpha - u i_X \alpha$ for an invariant form 
 $\alpha \in \Omega^*(W)^{S^1}$. 
 Now $d_{S^1}^2 = -u (d i_X + i_X d) = -u L_X = 0$ on 
 $\Omega^*(W)^{S^1} \otimes {\mathbb R}[u]$. 
 
 \begin{definition}
 The equivariant cohomology of the $S^1$-space $W$ is defined by
 $$H^*_{S^1}(W) = Ker d_{S^1} / Im d_{S^1}.$$
 \end{definition}
 
 From this definition, it is clear that $H^*_{S^1}(W)$ is a
 ${\mathbb R}[u]$-module. Furthermore, if $S^1$ acts trivially on $W$,
 then $H^*_{S^1}(W) \cong H^*(W) \otimes {\mathbb R}[u]$, a trivial module. 

 Let $P$ be
 a connected closed oriented manifold, and $\pi: E \rightarrow P$ be a 
 smooth complex vector bundle over $P$. Assume that there is an 
 $S^1$-action 
 on $E$ by bundle homomorphisms, which covers 
 an $S^1$-action on $P$. Following Atiyah-Bott \cite{Ati-Bot}, one can define
 the equivariant Euler class as
 $$\epsilon(E) = i^* i_* 1,$$
 where $i: P \rightarrow E$ is the zero section,  $i_*$ and $i^*$ are
 the push-forward and pullback homomorphisms in equivariant cohomology 
 respectively. It is routine to verify that
 $$\epsilon (E_1 \oplus E_2) = \epsilon(E_1) \epsilon(E_2)$$
 for two $S^1$ bundles $E_1$ and $E_2$ over $P$. We will be concerned with 
 the case when the action of $S^1$ on $P$ is trivial. In this case, by
 splitting principle \cite{Bot-Tu}, we can assume without loss of generality
 that $E$ has a decomposition as $S^1$ bundles 
 $$E = L_1 \oplus L_2 \oplus \cdots L_r,$$
 where each $L_j$ is a line bundle, such that 
 the action of $\exp(\sqrt{-1} t)$ on $L_j$ is multiplication by 
 $\exp (\sqrt{-1} m_j t)$, 
 for some weight $m_j \in {\mathbb Z}$. By formula $(8.8)$ in
 Atiyah-Bott \cite{Ati-Bot}, 
 $$\epsilon(L_j) = m_j u + c_1(L_j).$$
 Hence we have 
 $$\epsilon(E) = \prod_{j=1}^r (m_j u  + c_1(L_j)).$$

 \section{Localization formulas}
 
For an $S^1$-space $W$ with fixed point set $F$, let $\{ P_k \}$ be
the decomposition of $F$ into connected components. It is well-known that 
each $P_k$ is a smooth submanifold of $W$, hence $F$  has only finitely 
many components. The $S^1$-action on $W$ induces an
action on the normal bundle $\nu_k$ of $P_k$ in $W$. The equivariant Euler
class of $\nu_k$ can be computed as in $\S 1$. 
Now endow $W$ with an $S^1$-invariant metric.
Define a $1$-form $\theta$ on $W - F$ in the following way: $\theta(X) = 1$,
$\theta|_{X^{\perp}} = 0$. Here we use $X^{\perp}$ to denote the orthogonal
complement of $X$ in the tangent space. It is easy to see that
$\theta$ is a connection on the principal bundle $W - F \rightarrow 
(W - F) / S^1$. Following Kalkman \cite{Kal}, we define
for any $\alpha = \sum \alpha_j u^j \in \Omega^*(M)^{S^1} \otimes 
{\mathbb R}[u] $,
$$r(\alpha) = \sum \alpha_j (d \theta)^j - 
\theta \wedge (i_X \alpha_j) (d \theta)^j.$$
It is easy to see that $r(\alpha)$ is $S^1$-invariant and 
$i_X r(\alpha) = 0$. So $r(\alpha)$ is the lifting of a  form on  
$(W - F)/S^1$
via the projection $W - F \rightarrow (W - F) / S^1$.  Notice that there is an 
operator 
 $$\int_M: \Omega^*(W)^{S^1} \otimes
 {\mathbb R}[u] \rightarrow {\mathbb R}[u]$$
 induced by sending differential forms of degree $\dim (W)$ to its integral
 over $W$, and all other forms to zero.
 Now we can state a theorem due to Kalkman \cite{Kal} (which can be also 
obtained by Witten's localization principle \cite{Wit}).

\begin{theorem}
Let $W$ be an $S^1$-manifold with an invariant boundary $\partial W$, and  
fixed point set $F = \{ P_k \}$, such that $F \cap  \partial W = \emptyset$. 
Let $\alpha$ be an equivariant closed form on $M$ of total degree 
$dim (W) - 2$. Then
$$\int_{\partial W/ S^1} r(\alpha) = 
\sum_k \int_{P_k} \frac{\alpha u}{\epsilon (\nu_k)}.$$
\end{theorem}

We now give a construction of an equivariant closed form on $W$.  
Let $f: W \rightarrow {\mathbb R}$ be an $S^1$-invariant smooth function 
which vanishes near $F$, and $f \equiv 1$ outside a tubular neighborhood 
of $F$. Then $f\theta$ can be extended over $F$. It is straightforward  
to see that 
$d (f \theta)  - u ( -1 +i_X (f\theta)) = d (f \theta ) - u ( -1 + f)$
 is an equivariant closed form.
Assume now  that $\dim (W) = 2 (n + 1)$ and that $S^1$-action on the normal
bundles $\nu_k$ all have weight $1$. 
Let $\alpha = [ d (f \theta) - u ( -1 + f)]^n$, then 
near $\partial W$ we have $f \equiv 1$, and so

\[ r(\alpha) = r (( d \theta )^n )= (d\theta)^n.\]
Denote by $c$ the first 
Chern class of the principal $S^1$-bundle 
$\partial W \rightarrow \partial W /S^1$. Then in our normalization of 
$\theta$, $c = [- d \theta]$. 
An application of Theorem 2.1 then yields
$$\int_{\partial W/S^1} c^n = (-1)^n \sum_k \int_{P_k} 
\frac{u^{n+1}}{\sum_{j = 1}^{r_k}
c_{r_k - j}(\nu_k) u^j}, \eqno{(*)}$$
where $r_k$ is the complex rank of $\nu_k$. 

There is a slight generalization of the above formula. Let $\dim (W) = d + 2$,
$d$ is not necessarily even. Let $k$ be a number between $1$ and $d$,
which has the same parity as $d$. Let $\beta_1, \cdots, \beta_k$ be
$k$ $S^1$-invariant closed $1$-forms on $W$ such that for each 
$j = 1, \cdots, k$, $\beta_j|_{\partial W}$ is the pullback of a
$1$-form on $\partial W$, which we still denote by $\beta_j$. 
Now let $l = \frac{1}{2}(d - k)$, and let 
$$\alpha = \beta_1 \wedge \cdots \wedge \beta_k \wedge [d (f \theta) - u ( -1 + 
f)]^l.$$
Then near $\partial W$, 
$$r(\alpha) = r(\beta_1 \wedge \cdots \wedge \beta_k \wedge (d \theta)^l) 
=  \beta_1 \wedge \cdots \wedge \beta_k \wedge (d \theta)^l.$$
So by Theorem 2.1, we get 
$$\int_{\partial W/S^1} \beta_1  \wedge \cdots \wedge \beta_k \wedge c^l=  
\sum_k \int_{P_k} 
\frac{u^{l+1}\beta_1  \wedge \cdots \wedge \beta_k}{\sum_{j = 1}^{r_k}
c_{r_k - j}(\nu_k) u^j}. \eqno{(**)}$$

\section{Applications to Seiberg-Witten theory: A simple case}

Let $X$ be a closed oriented $4$-manifold. Given a Riemannian metric $g$ 
and a $Spin_c$ structure ${\mathcal S}$ on $X$, there are associated hermitian 
rank $2$ vector bundles $V_+$ and $V_-$, and a bundle isomorphism
$$\rho:  \Lambda_+ \rightarrow {\it su}(V_+),$$
where ${\it su}(V_+)$ is the bundle of anti-Hermitian traceless maps on 
$V_+$. 
The Seiberg-Witten equations are for a pair $(A, \Phi)$,
where $A$ is a unitary connection on $L = det (V_+)$, and $\Phi$ a 
section of $V_+$. For any fixed $\eta \in \Omega^+(X)$,
the perturbed Seiberg-Witten equations are
\begin{eqnarray}
\left\{
\begin{array}{rcl}
D_A \Phi & = & 0 \\
\rho(i F_A^+ + \eta) & = & (\Phi \otimes \Phi^*)_0
\end{array}
\right. \end{eqnarray}
These equations have a huge degree of  freedom. Let ${\mathcal A}$ denote the 
set of all unitary connections on $L$, ${\mathcal G}$ the group 
$Aut (L) = Map(X, S^1)$.
${\mathcal G}$ is called the gauge group. There is an action of ${\mathcal G}$ 
on
${\mathcal A} \times \Gamma(V_+)$, which preserves the Seiberg-Witten equations. 
 It is given by 
 $$g \cdot (A, \Phi) = (A - 2g^{-1} dg, g \Phi).$$
 This action is not free and  has two orbit types: if $\Phi \neq 0$, the 
 stabilizer of $(A, \Phi)$ is trivial; on the other hand, the stabilizer of
 $(A, 0)$ is $S^1$. To fix this problem, we choose an arbitrary point 
 $x_0 \in X$ and let ${\mathcal G}_0 = \{ g \in {\mathcal G} | g(x_0) = 1 \}$. 
 Then ${\mathcal G} = {\mathcal G}_0 \times S^1$, and furthermore, the action of 
 ${\mathcal G}_0$ on ${\mathcal A} \times \Gamma(V_+)$ is free. We then get the 
 residue $S^1$ action on the smooth 
 $({\mathcal A} \times \Gamma(V_+))/{\mathcal G}_0$.
  
Denote by $M({\mathcal S}, g, \eta)$ and $M^0({\mathcal S}, g, \eta)$ the 
quotients of the space ${\mathcal M}({\mathcal S}, g, \eta)$ of solutions 
to $(1)$ by $\mathcal G$ and ${\mathcal G}_0$
respectively.
They have the following well-known properties \cite{Kro-Mro}: 
\begin{itemize}
 \item{(a)} $M({\mathcal S}, g, \eta)$ and $M^0({\mathcal S}, g, \eta)$ are 
 compact in suitable
 topologies. 
 \item{(b)} For a generic choice of $(g, \eta)$, $M^0({\mathcal S}, g, \eta)$ 
 is a smooth manifold of
 dimension $d + 1 = 1 + \frac{1}{4}(c_1(L)^2 - 2 \chi(X) - 3 \tau(X))$, which 
can
 be oriented in a natural way. 
 \item{(c)} For a generic choice of $(g, \eta)$ with
 $2\pi c_1^+(L) \neq  \eta^h$, 
 the harmonic
part of $\eta$, $M^0({\mathcal S}, g, \eta)$ does not contain solutions with 
$\Phi = 0$ (called  reducible solutions).
 The $S^1$-action on $M^0({\mathcal S}, g, \eta)$ then gives rise
to a principal $S^1$-bundle $M^0({\mathcal S}, g, \eta) \rightarrow 
M({\mathcal S}, g, \eta)$. We call such a choice of $(g, \eta)$ a good choice.
When $b_2^+(X) > 0$, there are good choices. 
\item{(d)} For two good choices
$(g_0, \eta_0)$ and $(g_1, \eta_1)$, there is a path $(g_t, \eta_t)$ joining
them, such that $M^0({\mathcal S}, g_t, \eta_t)$, $0 \leq t \leq 1$, form an
oriented cobordism $W$ between $M^0({\mathcal S}, g_0, \eta_0)$ and
$M^0({\mathcal S}, g_1, \eta_1)$. When $b_2^+(X) > 1$, it is possible to choose
the path such  that none of $M^0({\mathcal S}, g_t, \eta_t)$ admits a reducible 
solution. 
\end{itemize}

For a good choice $(g, \eta)$, the Seiberg-Witten invariant is defined 
as follows: 
(a) if $d < 0$, $SW({\mathcal S}, g, \eta) = 0$; 
(b) if $d = 0$, then $M({\mathcal S}, g, \eta)$ is a finite union of signed 
points, and $SW({\mathcal S}, g, \eta)$ is the sum of the corresponding $\pm 
1$'s; 
(c) if $d > 0$, the Seiberg-Witten invariant can be defined as the coupling
of the fundamental class of $M({\mathcal S}, g, \eta)$ with the suitable power 
of
the first Chern class of the principal $S^1$-bundle $M^0({\mathcal S}, g, \eta)
\rightarrow M({\mathcal S}, g, \eta)$.
So when $b_2^+(X) > 1$, $SW({\mathcal S}, g, \eta)$ does not depends on the 
good choice $(g, \eta)$ and is then a diffeomorphism invariant. 
However, if $b_2^+(X)  = 1$, $c_1^+(L) = \eta^h$ defines a hypersurface
in the space of $(g, \eta)$'s. It is called a ``wall", since it divides
the space of $(g, \eta)$'s into two connected components, called ``chambers".
For two good choices $(g_0, \eta_0)$ and $(g_1, \eta_1)$, the Seiberg-Witten
invariants are the same. However, when they lie in different chambers the 
invariants
may differ. A formula relating the Seiberg-Witten invariants for 
good choices in different chambers is called a wall crossing formula. 

Since the invariant is nontrivial only if $d$ is odd, the above formula 
for $d$ shows that the only interesting case is when $b_2^+(X) = 1$, and 
$b_1(X)$ is even. The wall crossing formula of Seiberg-Witten invariants
in the case $b_2^+ = 1$, $b_1 = 0$ and $d = 0$ was obtained by
Witten \cite{Wit} and Kronheimer-Mrowka \cite{Kro-Mro}. The general wall
crossing formula, proved by Li-Liu \cite{Li-Liu}, can be stated as  follows

\begin{theorem}
Let $X$ be a closed oriented $4$-manifold with $b_2^+ = 1$ and $b_1$ even,
${\mathcal S}$ a $Spin_c$ structure with $det(V_+) = L$, such that
$c_1(L)^2 - (2 \chi(X) + 3 \tau(X)) \geq 0$, then for any two good choices 
$(g_0,\eta_0)$ and $(g_1, \eta_1)$ in two different chambers, the Seiberg-Witten
invariants $SW({\mathcal S}, g_0, \eta_0)$ and $SW({\mathcal S}, g_1, \eta_1)$ 
differ
by
$$\pm \int_{T^{b_1}}
(\frac{1}{4} \Omega^2 \cdot c_1(L)[X])^{b_1/2}/(b_1/2)!,$$
where 
$$\Omega = c_1({\mathcal U}) = \sum_i x_i \cdot y_i,$$
and ${\mathcal U}$ is the universal flat line bundle over $T^{b_1} \times 
M$, $\{y_i \}$ is any basis of $H^1(X; {\mathbb Z})$ modulo torsion, and 
$\{ x_i \}$ is the dual basis in $H^1(T^{b_1}; {\mathbb Z})$.
\end{theorem}

We will now reprove this theorem by the method described in $\S 2$. 
Take a path $(g_t, \eta_t)$ that goes through the wall transversally once.
Then the $S^1$-action on the induced cobordism $W$ has only one component 
in the fixed point set $F$, namely the set of reducible solutions, which are 
parameterized by the torus $T^{b_1} = H^1(X; {\mathbb R})/H^1(X; {\mathbb Z})$.
We shall assume that for each reducible solution $(A, 0)$, $Coker D_A = 0$. 
(The general case can be modified by the method of 
Li-Liu \cite{Li-Liu}, p. 808.) Under this assumption, the normal bundle of 
$F$ in $W$ is given by the index bundle $ind $, whose fiber at each  
$(A, 0) \in F$ is given by $Ker D_A$ (cf. Li-Liu \cite{Li-Liu}.)  
It is clear that the $S^1$-action on this normal bundle has only weight 
$1$. Now we  use formula $(*)$ in $\S 2$  to get 
 \begin{eqnarray}
 & & SW({\mathcal S}, g_1, \eta_1) - SW({\mathcal S}, g_0, \eta_0) 
 = \int_{\partial W/S^1} c^{d/2} \\
 & = & \pm \int_{T^{b_1}} \frac{u^{(d+1)/2}}{\sum_{j=1}^r c_{r-j}(ind) u^j} 
 \end{eqnarray}
 where $r$ is the complex rank of $ind$, so $2r + b_1 = d+2$. In the proof 
 of Lemma 2.5 in \cite{Li-Liu}, Li-Liu derived, by
 Atiyah-Singer family index theorem, that
 \begin{eqnarray*}
c_1(ind) & = & \frac{1}{4}\Omega^2 \cdot c_1(L)[X], \\
c_j(ind) & = & \frac{1}{j!}c_1(ind)^j.
\end{eqnarray*}
Plugging the above equalities into (4),  we see that the difference between the 
two Seiberg-Witten invariants
is 
\begin{eqnarray*}
& & \pm \int_{T^{b_1}} \frac{u^{(d+2)/2}}{u^r exp (c_1(ind)/u)} \\
& = & \pm \int_{T^{b_1}} u^{-r+(d+2)/2} exp ( - c_1(ind)/u) \\
& = & \pm \int_{T^{b_1}} u^{-r+(d+)/2 - b_1/2} c_1(ind)^{b_1/2}/(b_1/2)! \\
& = & \pm \int_{T^{b_1}} c_1(ind)^{b_1/2}/(b_1/2)!
\end{eqnarray*}
This completes the proof of Theorem 3.1.

\section{Applications to Seiberg-Witten theory: The general case}

Okonek-Teleman \cite{Oko-Tel} extended the definition of Seiberg-Witten 
invariants. They also proved a
wall crossing formula for such general Seiberg-Witten invariants. 
In this section, we will give an equivalent definition of the
general Seiberg-Witten invariants, which is along the line of 
the discussions in the preceding sections. We then reprove Okonek-Teleman's
formula by the localization formula $(**)$.

We use the notations of $\S 3$. Let $L \rightarrow X$ be the Hermitian line
bundle associated to a fixed $Spin_c$-structure ${\mathcal S}$. 
For a good choice 
$(g, \eta)$, let $\pi_2: {\mathcal M}({\mathcal S}, g, \eta) 
\times X \rightarrow X$ be the projection onto the second factor. Consider
the pullback line bundle
$\pi^*_2L$. The group ${\mathcal G}$ acts freely on $\pi_2^* L$, which 
covers a free action of ${\mathcal G}$ on ${\mathcal M}({\mathcal S}, g, \eta) \times
X$. There is therefore a quotient line bundle 
$${\mathcal L} \rightarrow M({\mathcal S}, g, \eta) \times X.$$
We now define a group homomorphism
$\mu: H_1(X; {\mathbb Z}) / Tor \rightarrow H^1(M({\mathcal S}, g, \eta); 
{\mathbb R})$
by 
$$\mu([A]) = \int_A c_1({\mathcal L}),$$
where $A$ is a loop in $X$, and $[A]$ its homology class. It is easy to
see that this is well-defined. 

 Let $d = \frac{1}{4}(c_1(L)^2 - 2 \chi (X) - 3 \tau (X))$. 
 When $ d < 0$, the Seiberg-Witten
invariant $SW({\mathcal S}, g, \eta)$ is defined to be zero. When $d = 0$,
it is defined as in $\S 3$. When $d > 0$, $SW({\mathcal S}, g, \eta)$ is
defined as a linear map
$$\Lambda^* (H_1(M, {\mathbb Z})/ Tor) \rightarrow {\mathbb R}.$$
More precisely, for $ 0 \leq k \leq \min \{b_1, d\}$, and $k$ has the same
parity as $d$,
$$SW({\mathcal S}, g, \eta)([A_1], \cdots, [A_k]) = \int_{M({\mathcal S}, g, 
\eta)}
\mu([A_1]) \wedge \cdots \wedge \mu([A_k]) \wedge c^l,$$
where $l = \frac{1}{2}(d - k)$ and $c$ is as in $\S 3$.
For all other values of $k$, the invariant is defined to be zero. 
  We remark that
these invariants are actually integer-valued, even though we define them as 
integrals of differential forms.

Okonek-Teleman's wall crossing formula can be stated as the following

\begin{theorem} For a fixed $Spin_c$-structure ${\mathcal S}$ on a connected 
closed oriented $4$-manifold $X$ with $b_2^+ = 1$, $d = \frac{1}{4} (c_1(L)^2
- 2 \chi(X) - 3 \tau(X)) \geq 0$, the Seiberg-Witten invariants for
two good choices $(g_0, \eta_0)$ and $(g_1, \eta_1)$ in two different
chambers are related by
\begin{eqnarray*}
& & SW({\mathcal S}, g_1, \eta_1) ([A_1], \cdots, [A_k]) - 
    SW({\mathcal S}, g_0, \eta_0) ([A_1], \cdots, [A_k]) \\
& = & \pm \int_{T^{b_1}}\mu([A_1]) \wedge \cdots \wedge \mu([A_k]) \wedge 
(\frac{1}{4} \Omega^2 \cdot c_1(L) [X])^{(b_1 - k)/2} / ((b_1 - k)/ 2)!
\end{eqnarray*}
where $[A_1], \cdots, [A_k]$ are homology classes in $H_1(X; {\mathbb Z})/Tor$,
$\Omega$ is as in Theorem 3.1, $ 0 \leq k \leq \min \{ b_1, d \}$,
and $k$ has the same parity as $b_1$. 
\end{theorem}

The proof is similar to the proof of Theorem 3.1. To start with, we give a
parameterized version of the construction for ${\mathcal L}$. Take a
path $(g_t, \eta_t)$ as in $\S 3$. Consider the
infinite dimensional cobordism
$${\mathcal W} = \bigcup_{t} {\mathcal M}({\mathcal S}, g_t, \eta_t).$$
Consider the pullback line bundle $\pi_2^*L$ on ${\mathcal W} \times X$,
where $\pi_2$ is again the projection onto the second factor.
Modulo the action by ${\mathcal G}_0$, we get a quotient line bundle
$${\mathcal L}_0 \rightarrow W \times X.$$
This is actually an $S^1$-bundle, since there is  the residue action by
$S^1 = {\mathcal G} / {\mathcal G}_0$. Since this $S^1$-action on $\partial W$
is free, it is straightforward to see that when restricted to 
$\partial W = M_0({\mathcal S}, g_0, \eta_0) \sqcup M_0({\mathcal S}, g_1, 
\eta_1)$,
${\mathcal L}_0$ can be identified with the pullback of $\mathcal L$ on 
$M({\mathcal S}, g_j, \eta_j)$ via the projection $M_0({\mathcal S}, g_j, 
\eta_j)
\rightarrow M({\mathcal S}, g_j, \eta_j)$, for $j = 1, 2$. 
Endowing ${\mathcal L}_0$ an $S^1$-invariant unitary connection, we then see 
that
$c_1({\mathcal L}_0)$ is represented by an $S^1$-invariant closed $2$-form. 
Define $\mu_0: H_1(X; {\mathbb Z}) / Tor \rightarrow H^1(W; {\mathbb R})$
by 
$$\mu_0([A]) = \int_A c_1({\mathcal L}_0).$$
It is easy to see that when restricted to $\partial W$, 
$\mu_0([A  ])$ is the pullback of $\mu([A])$ on $M({\mathcal S}, g_0, \eta_0)$
and $M({\mathcal S}, g_1, \eta_1)$ respectively.
Now let $\beta_j = \mu_0([A_j])$, for $j = 1, \cdots, k$, a computation similar
to the one in $\S 3$ by formula $(**)$ then proves Theorem 4.1.

{\em Acknowledgement}. The effort of applying equivariant cohomology 
to wall crossing was initiated in the second author's joint effort
with Rugang Ye in understanding the contribution of reducible 
solution in $3$-dimensional Seiberg-Witten theory while he was visiting 
University of California at Santa Barbara.
The present work is carried out during his visit at Texas A\& M
University. He likes to thank Mathematics Department at both 
universities for
hospitality and financial support. He also thanks Xianzhe Dai, 
Doug Moore, Guofang Wei and Rugang Ye for helpful discussions and
encouragement, and Blaine Lawson for teaching him equivariant cohomology.

\bigskip

\noindent Department of Mathematics,
Texas A\&M University,
College Station, TX 77843

\noindent {\it E-mail address:}\ cao@math.tamu.edu,  zhou@math.tamu.edu


\begin{thebibliography}{99}


\bibitem{Ati-Bot} M.F. Atiyah, R. Bott, {The moment map and equivariant
cohomology}, {\bf Topology 23} (1984), 1-28.

\bibitem{Ber-Get-Ver} N. Berline, E. Getzler, M. Vergne, {\bf 
Heat kernels and Dirac operator}, Springer, 1992.

\bibitem{Bot-Tu} R. Bott, L. Tu, {\bf Differential forms in algebraic topology,
GTM 82}, Springer, 1982.

\bibitem{Got} L. G\"{o}ttsche, {\em Modular forms and Donaldson invariants 
for $4$-manifolds with $b\sb {+}=1$}, {\bf J. Amer.
Math. Soc. 9} (1996), no. 3, 827--843. 

\bibitem{Kal} J. Kalkman, {\em Cohomology rings of symplectic quotients},
{\bf J. Reine Angw. Math. 458} (1995), 37-52.

\bibitem{Kro-Mro} P. Kronheimer, T. Mrowka, {\em The genus of embedded surfaces
in the projective plane}, {\bf Math. Res. Letters 1} (1994), 797-808.

\bibitem{Li-Liu} T.J. Li, A. Liu, {\em General wall crossing formula},
{\bf Math. Res. Letters 2}, (1995), 797-810.

\bibitem{Oko-Tel} C. Okonek, A. Teleman, {\em Seiberg-Witten 
invariants for manifolds with $b\sb +=1$ and the
universal wall crossing formula}, {\bf Internat. J. Math. 7 }(1996), 
no. 6, 811--832.

\bibitem{Qui} V. Mathai, D. Quillen, {\em Superconnections, Thom classes
and equivariant differential forms}, {\bf Topology 25} (1986), 85-110.

\bibitem{Li-Tia} G. Tian, private communication.

\bibitem{Wit} E. Witten, {\em Monopoles and $4$-manifolds}, {\bf Math. 
Res. Letters 1} (1994), 769-796.
\end{thebibliography}
\end{document}